# SHARP ESTIMATES FOR SPECTRUM OF CAUCHY OPERATOR AND PERIODS OF SOLUTIONS OF LIPSCHITZ DIFFERENTIAL EQUATIONS

A.A. Zevin

*Institute Transmag, Academy of Sciences of Ukraine*
*Piesarzhevsky 5. 49005, Dnepropetrovsk, Ukraine*

Phone: (38056)732-30-55
Fax: (38056)370-21-86
e-mail: zevin@westa-inter.com

**Abstract**   Estimates for the spectrum of the Cauchy operator and logarithms of solutions of non-autonomous differential equations in the space $\mathbb{C}^n$, expressed in an arbitrary matrix norm, are found. For equations with periodic coefficients, the lower bound for the periods of "oscillatory" solutions is obtained. Similar results are derived for a nonlinear periodic equation with zero equilibrium position; in the autonomous case, they are valid for any periodic solution. All the estimates are attained at a scalar equation and, thereby, are accurate for any norm.  The estimates are extended to equations with derivatives of any order.

Using the obtained results, a condition for the uniqueness of an oscillatory solution is found; stability criteria for a linear Hamiltonian system with periodic coefficients and for periodic solutions of a nonlinear system (expressed through the period and any norm of the Hessian) are obtained. Types of matrices and functions, for which the Euclidean norm has a minimum value (and, therefore, provides the most accurate estimates for a specific equation), are established.

## 1. Introduction

We consider the following differential equations

$$\dot{z} = A(t)z, \qquad (1.1)$$

$$\dot{z} = f(z,t), \tag{1.2}$$

$$f(z,t) = f(z,t+T), \quad f(0,t) = 0$$

and

$$\dot{z} = f(z). \tag{1.3}$$

Here $z \in \mathrm{C}^n$, the matrix $A(t)$ is piecewise continuous, the functions $f(z)$ and $f(z,t)$ satisfy the Lipschitz conditions,

$$\|f(z') - f(z'')\| \le L\|z' - z''\| \tag{1.4}$$

and

$$\|f(z',t) - f(z'',t)\| \le L(t)\|z' - z''\| \tag{1.5}$$

where $\|\cdot\|$ is any norm in $\mathrm{C}^n$.

First, the localization of the spectrum, $\rho_k(t)$, $k = 1,...,n$, of the Cauchy operator of equation (1.1) is considered. The known results of this kind are set for the Hilbert space (Daletskii and Krein, [1]). Below sharp estimates for $|\ln \rho_k(t)|$, expressed in

$$N(t) = \int_0^t \|A(s)\| ds, \tag{1.6}$$

are obtained.

Further, we find the upper bound (2.7) for $\|\ln z(t)\|$ where $z(t)$ is a solution, of equation (1.1). In contrast to the well-known inequality,

$$\|z(t)\| \le \|z(0)\| \exp N(t), \tag{1.7}$$

it provides bound not only for the modules, but also for the angular components, $\varphi_k(t) = \arg z_k(t)$, of the solution $z(t)$.

For the case $A(t) = A(t+T)$, we seek for the minimum value of $T$, for which there exists a $T$-periodic solution (apparently, such a problem is posed here for the first time). The set of such solutions is classified as follows.

**Definition**. A solution $z(t) = z(t+T)$ is called oscillatory if for some $k$ and $t_1, t_2 \in [0,T)$,

$$|\varphi_k(t_2) - \varphi_k(t_1)| = \pi. \tag{1.8}$$

Otherwise, the solution $z(t)$ is referred to as non-oscillatory.

Condition (1.8) means that the component $z_k(t)$ does not lie entirely in some half-plane. In particular, condition (1.8) holds when $z_k(t)$ has a zero mean value on the interval $[0,T]$. The last, for example, is true for odd or anti-periodic ($z(t) = -z(t+T/2)$) solutions.

We show that the periods of oscillatory (unlike, non-oscillatory) solutions are bounded from below; the corresponding sharp bound, $T_*$, is determined from the relation $N(T) = 2\pi$. This result is extended to equation (1.2) (it is only necessary to replace $\|A(t)\|$ by $L(t)$).

For equation (1.3), the problem is to find the sharp lower bound, $T_*(L)$, for periods of solutions $z(t) = z(t+T) \neq const$ (since $T_*(L) = T_*(1)/L$, it is sufficient to find $T_*(1)$). Such problem was first posed by Yorke [2], who found that for equation (1.3) in the space $\mathrm{R}^n$ with the Euclidean norm,

$$T_*(1) = 2\pi. \tag{1.9}$$

The same result holds in the Hilbert space (Lasota and Yorke, [3]) as well as for the equation

$$z^{(r)} = f(z) \tag{1.10}$$

with any $r$ (Mawhin and Walter, [4]).

The proofs of these results essentially use the presence of a scalar product in the definition of the norm, so they are not applicable to other metrics. Using a different approach, it is proved (Zevin, [5]) that for equation (1.10) with even $r$ in the space $\mathrm{R}^n$, equality (1.9) is valid for a wide class of norms, including the most common ones. The same result was claimed for odd $r$ [5, 6], however, the corresponding proofs contain errors (the author is very grateful to M.Nieuwenhuis and J. Robinson which pointed out them). The refined proof for equation (1.10) in the space $\mathrm{C}^n$ with any $r$ and arbitrary norm is given below.

Note that for equation (1.3) in the general Banach space, $T_*(1) = 6$ (Busenberg, Fisher and Martelli, [7]).

The obtained results are extended to equations

$$z^{(r)} = A(t)z \tag{1.11}$$

and

$$z^{(r)} = f(z,t). \tag{1.12}$$

In particular, it is found that for $N(T)T^{-1} = L = 1$, the minimal periods of oscillatory solutions to equations (1.11), (1.12) and non-constant solutions to equation (1.10) are equal $2\pi$, regardless of $r$, $n$ and the accepted norm.

All the estimates hold true in the space $\mathrm{R}^n$ and, in the case of even $r$, remain sharp; for odd $r$, they are sharp for universally adopted norms.

As an application, a condition for the uniqueness of an oscillatory solution to system (1.12) with an additional perturbation $p(t) = p(t+T)$ is obtained. The known criteria by Krein [1] for stability of a Hamiltonian system with periodic coefficients in the Hilbert space are generalized to the space $\mathrm{C}^n$ with any norm; for a nonlinear Hamiltonian system, a condition for the stability of a periodic solution is found.

For a specific equation, the obtained estimates depend on the choice of a metric; in this connection, classes of matrices and functions, for which the Euclidean metric provides the best estimates, are indicated.

## 2. Main results

Let $W(t,0)$ and $\rho_k(t)$, $k=1,\ldots n$ be the Cauchy operator of equation (1.1) and its eigenvalues. Since $W(0,0) = I$ is the identity matrix, then $\rho_k(0) = 1$; having put $\arg \rho_k(0) = \varphi_k(0) = 0$, we shall determine $\varphi_k(t)$ by continuity.

**Theorem 1**. In system (1.1),
$$|\ln \rho_k(t)| \le N(t), \; k=1,\ldots n. \tag{2.1}$$

**Proof.** Set $W(t,0) = \exp K(t)$, then $\ln \rho_k(t)$ are the eigenvalues of the matrix $K(t)$. Hence,
$$|\ln \rho_k(t)| \le \|K(t)\|. \tag{2.2}$$

Let $t_k < t_{k+1}$, $t_0 = 0$, $t_{p+1} = t$, then
$$W(t,0) = \prod_{k=0}^{p} W(t_{k+1}, t_k).$$

Approximate $A(t)$ with a piecewise constant function $A_*(t) = A(t_k)$ for $t \in [t_k, t_{k+1})$. For the corresponding Cauchy operator, $W_*(t,0)$, one has
$$\|W_*(t_{k+1}, t_k)\| = \|\exp[A(t_k)(t_{k+1} - t_k)]\| \le \exp(\|A(t_k)\|(t_{k+1} - t_k)),$$

so

$$\|W_*(t,0)\| \leq \exp\left(\sum_{k=0}^{p}\|A(t_k)\|(t_{k+1}-t_k)\right).$$

Put $k \to \infty$, $|t_{k+1}-t_k| \to 0$, then $\lim W_*(t,0) = W(t,0)$. Hence,

$$\|W(t,0)\| \leq \exp\int_0^t \|A(s)\|ds = \exp N(t)$$

and, therefore,

$$\|K(t)\| \leq N(t). \tag{2.3}$$

Inequality (2.1) follows from (2.2) and (2.3). □

In particular, from (2.1) one has

$$|\arg \rho_k(t))| \leq N(t) \tag{2.4}$$

(for equation (1.1) in the Hilbert space, such inequality was established in [1]).

For any norm, the equality in (2.1) is achieved at the system of uncoupled equations

$$\dot{z}_k = \exp(is_k)\|A(t)\|z_k, \quad s_k \in [0,1], \quad k=1,...,n \tag{2.5}$$

The following theorem provides an upper bound for $\|\ln z(t))\|$ ($\varphi_k(0)=0$, $k=1,...n$).

**Theorem 2**. In system (1.1),

$$\|\ln z(t)\| \leq \ln\|z(0)\| + N(t) \tag{2.6}$$

**Proof**. Putting $z(t) = \exp[\ln z(t)]$, we rewrite the relation $z(t) = W(t,0)z(0)$ in the form

$$\exp[\ln z(t)] = \exp K(t) \exp[\ln z(0)].$$

Therefore,

$$\exp\|\ln z(t)\| \leq \exp\|K(t)\|\exp\|\ln z(0)\| = \exp\|K(t)\|\exp\ln\|z(0)\|,$$

whence one has

$$\|\ln z(t)\| \leq \|K(t)\| + \ln\|z(0)\| \leq N(t) + \ln\|z(0)\|.$$

□

The equality in (2.6) is also achieved at system (2.6).

To compare Inequalities (2.6) and (1.7), we rewrite the latter in the form

$$\ln\|z(t)\| \leq \ln\|z(0)\| + N(t), \tag{2.7}$$

Since

$$|\ln z_k(t)| = |\ln|z_k(t)| + i\varphi_k(t)| \geq |\ln|z_k(t)||,$$

then inequality (2.6) is more precise.

Let us evaluate $\|\varphi(t)\|$. Without loss of generality, assume $\|z(0)\|=1$, then $\ln\|z(0)\|=0$ and (2.6) implies

$$\|\varphi(t)\| \leq \|\ln z(t)\| \leq N(t). \tag{2.8}$$

Suppose now that

$$A(t) = A(t+T) \tag{2.9}$$

The following theorem provides a lower bound for the periods of oscillatory solutions of equation (1.1).

**Theorem 3**. If

$$N(T) < 2\pi, \tag{2.10}$$

then system (1.1), (2.11) has no $T$ - periodic oscillatory solutions.

**Proof**. Let $z(t) = z(t+T)$ be an oscillatory solution. Put in (1.8) $t_2 = 0$, then $|\varphi_k(t_1)| = \pi$ and $|\varphi_k(t_1+T)| = \pi + 2\pi p$ where $p$ is an integer. From (2.8) we obtain the inequality

$$2\pi \leq |\varphi_k(t_1)| + |\varphi_k(t_1+T)| \leq \|\varphi(t_1)\| + \|\varphi(t_1+T)\| \leq \int_{t_1}^{0}\|A(t)\|dt + \int_{0}^{t_1+T}\|A(t)\|dt = N(T)$$

which contradicts (2.10). □

The lower bound for the period of oscillatory solutions is determined from the relation

$$N(T) = 2\pi. \tag{2.11}$$

As is clear from the proof, equality (2.11) is valid if only

$$\left|\frac{d \arg z_k(t)}{dt}\right| \equiv \left|\frac{d \ln z_k(t)}{dt}\right| \equiv \|A(t)\|. \tag{2.12}$$

For any norm, identity (2.12) holds in the system

$$\dot{z}_k = \pm i\|A(t)\|z_k \qquad k=1,...,n \tag{2.13}$$

Note that for the periods of non-oscillatory solutions, a similar bound does not exist. E.g., the real system

$$\dot{x}_1 = x_2, \quad \dot{x}_2 = -a(t)x_1, \tag{2.14}$$

$$a(t) = a(t+T), \int_{0}^{T} a(t)dt < 0.$$

may have a $T$ - periodic solution with an arbitrarily small $T$ [8]. Herewith, the component $x_1(t)$ is sign-constant, so the solution $z_1(t) = x_1(t) + ix_2(t)$ of the respective complex equation is non-oscillatory.

Consider now equation (1.2). Put

$$N_L(t) = \int_0^t L(s)ds. \qquad (2.15)$$

**Theorem 4**. If

$$N_L(T) < 2\pi, \qquad (2.16)$$

then equation (1.2) has no $T$ - periodic oscillatory solutions.

**Proof.** Using integral theorem on finite increments and taking into account that $f(0,t) = 0$, we obtain

$$f(z,t) = f(z,t) - f(0,t) = C(z,t)z, \quad C(z,t) = \int_0^1 f_z(sz,t)ds. \qquad (2.17)$$

Therefore, a solution $z(t)$ of equation (1.2) satisfies also

$$\dot{z} = C(t)z. \qquad (2.18)$$

As is known, the function $L(t)$ in (1.5) can be defined as

$$L(t) = \sup_{z \in C^n} \|f_z(z,t)\|. \qquad (2.19)$$

Obviously, $C(t) \leq L(t)$, so Theorem 4 follows from Theorem 3. □

The lower bound for the periods of the oscillatory solutions is defined here by

$$N_L(T) = 2\pi. \qquad (2.20)$$

The following theorem gives the sharp estimate for the period of non-constant solutions of equation (1.3).

**Theorem 5**. In system (1.3), (1.4) with any norm, the period of a solution $z(t) = z(t+T) \neq const$,

$$T \geq \frac{2\pi}{L}. \qquad (2.21)$$

**Proof.** Since equation (1.3) is autonomous, the function $y(t) = \dot{z}(t)$ satisfies

$$\dot{y} = A(t)y, \qquad (2.22)$$

where $A(t) = f_z(z(t))$.

The Lipschitz constant of the function $f(z)$ can be defined as

$$L = \sup_{z \in C^n} \|f_z(z)\|, \qquad (2.23)$$

so

$$\|A(t)\| \leq L, \quad N(T) \leq LT. \qquad (2.24)$$

Since

$$\int_0^Y y(t)dt = z(T) - z(0) = 0,$$

the solution $y(t)$ is oscillatory. Hence, inequality (2.21) follows from (2.24) and Theorem 3. □

The equality in (2.21) is achieved at the system

$$\dot{z}_k = \pm iLz_k, \qquad k = 1,...,n. \qquad (2.25)$$

The next theorem generalizes Theorem 3 to equation (1.11).

**Theorem 6**. If for some $c > 0$,

$$N(T,c) = \int_0^T \max[c, c^{1-r}\|A(t)\|]dt < 2\pi, \qquad (2.26)$$

then equation (1.11) has no $T$ - periodic oscillatory solutions.

**Proof**. Putting

$$v_1 = z, \quad \dot{v}_{k-1} = cv_k, \quad k = 2,...,r-1, \quad \dot{v}_r = c^{1-r}A(t)v_1, \quad v = [v_1,...,v_r],$$

we obtain the system of order $n \times r$,

$$\dot{v} = V(t.c)v. \qquad (2.27)$$

Obviously,

$$\frac{\|V(t,c)v\|}{\|v\|} = \frac{c\sum_{k=2}^{r}\|v_k\| + c^{1-r}\|A(t)v_1\|}{\sum_{k=1}^{r}\|v_k\|}$$

whence we find

$$\|V(t,c)\| = \max[c, c^{1-r}\|A(t)\|] . \qquad (2.28$$

Application of Theorem 3 to system (2.27), (2.26) proves the theorem. □

The equality

$$N(T,c) = 2\pi \qquad (2.29)$$

gives a lower bound, $T(c)$, for periods of oscillatory solutions. It can be refined by calculating

$$T_* = \max_c T(c). \qquad (2.30)$$

If $\|A(t)\| = \|A\|$ is constant, the value $T_*$ is attained for $c = \|A\|^{1/r}$; wherein

$$\|V(t,c)\| = \|A\|^{1/r}, \quad T_* = \frac{2\pi}{\|A\|^{1/r}} . \qquad (2.31)$$

This bound is accurate: it is achieved at the equation

$$z^{(r)} = \pm i^r \|A\| z. \qquad (2.32)$$

If $\|A(t)\|$ is non-constant and $r > 1$, estimate (2.30) is not reached (for $r = 1$, $\|V(t,c)\| = \|A(t)\|$, so the exact bound is determined from the relation (2.11)).

Analogously to (2.18), a solution $z(t)$ of equation (1.12) satisfies also the equation

$$\dot{z}^{(r)} = C(t)z, \qquad (2.33)$$

where $C(t) \leq L(t)$. This implies the following theorem.

**Theorem 7**. If for some $c > 0$,

$$N(T,c) = \int_0^T \max[c, c^{1-r} L(t)] dt < 2\pi, \qquad (2.34)$$

then equation (1.12) has no $T$ - periodic oscillatory solutions.

For equation (1.10), the solution $z(t) = z(t+T) \neq const$ meets the variational equation

$$y^{(r)} = A(t)y \qquad (2.35)$$

where $\|A(t)\| \leq L$. Since (1.10) is autonomous, the solution $y(t) = \dot{z}(t)$ has zero mean value and, hence, is oscillatory. So, the greatest lower bound for $T$,

$$T_* = \frac{2\pi}{L^{1/r}}. \qquad (2.36)$$

It is achieved at the system

$$z_p^{(r)} = \pm i^r L z_p, \qquad p = 1,\ldots,n. \qquad (2.37)$$

## 3. Discussion

If $z \in \mathbb{R}^n$, equations (1.1)-(1.3) may be written in the complex form. (in the case of odd $n$, one may add the equation $\dot{z}_{n+1} = 0$). Hence, the above estimates hold for real systems; however, their accuracy for a specific norm needs to be checked as follows.

As is shown above, in the case $z \in \mathbb{C}^n$, the sharp estimates are achieved at a system of the kind (2.37). If $r = 2q$, $q = 1,2,\ldots$, it has the real solution

$$z_p(t) = a_p \sin(L^{1/r} t) + b_p \cos(L^{1/r} t) \qquad (3.1)$$

with the period $T = 2\pi / L^{1/r}$. Consequently, for even $r$, the obtained estimates remain sharp in the space $\mathrm{R}^n$ for any norm.

In the case $r = 2q - 1$, (2.37) is equivalent to the real system

$$x_{q1}^{(r)} = (-1)^q L x_{q2}, \quad x_{q2}^{(r)} = -(-1)^q L x_{q1}, \quad q = 1,2,... \tag{3.2}$$

which also has a solution with the period $T = 2\pi / L^{1/r}$. Therefore, the obtained estimates remain sharp in $\mathrm{R}^n$, if for some $q = p$, the norm of the right side of (3.2) equals $L$, i.e.,

$$\|(x_{p2}, -x_{p1})\| = \|(x_{p1}, x_{p2})\|. \tag{3.3}$$

Really, setting $x_{q1}^{(r)} = 0$, $x_{q2}^{(r)} = 0$ for $q \neq p$, we obtain the system with the period $T = 2\pi / L^{1/r}$ and the Lipschitz constant $L$. .

Note that the most common norms, in particular,

$$\|x\|_E = \left(\sum_{p=1}^n x_p^2\right)^{1/2}, \quad \|x\|_S = \max_p |x_p|, \quad \|x\|_O = \sum_{p=1}^n |x_p|, \tag{3.4}$$

satisfy condition (3.3).

For a specific equation, the obtained results provide the upper bounds for $|\ln \rho_k(t)|$, $\|\ln z(t))\|$ and lower bounds for the period $T$ of oscillatory solutions. Clearly, in all the cases, the most accurate bound is achieved at the metric, for which the corresponding value $N$, $N_L$ or $L$ is minimal.

In [6] it is shown that for Hermitian (anti-Hermitian) matrices, such values are achieved in the Euclidean norm. The following lemma generalizes this result to the normal matrices ($AA^* = A^*A$)..

Let $\|A\|_E$ and $\|A\|$ be the Euclidean and any other norm of a matrix $A$.

**Lemma 1**. For a normal matrix,

$$\|A\|_E \leq \|A\|. \tag{3.5}$$

**Proof**. By definition, $\|A\|_E = \lambda_*^{1/2}$, where $\lambda_*$ is the largest eigenvalue of the matrix $A^*A$. Let $\mu_i$ be the eigenvalues of $A$. As is known [9], a normal matrix can be represented in the form $A = U^* \Delta U$, where $U$ is a unitary matrix and $\Delta = \mathrm{diag}[\mu_i]$, Hence, $\|A\|_E = \max|\mu_i| \leq \|A\|$ what proves (3.5). □

It is easy to show that Lemma 1 cannot be extended to a general matrix.

Using this lemma, let us indicate a class of functions for which the Euclidean Lipschitz constant has the smallest value. Let

$$f(x,t) = V_x(x,t), \quad x \in \mathrm{R}^n \tag{3.6}$$

where the function $V(x,t)$ is twice differentiable in $x$. Since the matrix $f_x(x,t) = V_{xx}(x,t)$ is symmetric, its Euclidean norm has the minimal value. Taking into account representation (2.23), we find that the Euclidean Lipschitz constant of the function $f(x,t)$ is also minimal.

The same result is true for the function

$$f(x,t) = RV_x(x,t), \quad x \in \mathrm{R}^n \tag{3.7}$$

where $R$ is an orthogonal matrix. Indeed, it is easy to show that the matrix $f_x(x,t)$ is normal.

Thus, for systems with the specified functions and matrices, the best estimates are provided with the Euclidean norm.

Note that the last conclusion is applicable to other problems for which some estimate is valid for any matrix norm or Lipschitz constant (e.g., such estimates are known for the upper and lower Lyapunv exponents of equation (1.1) [1]).

## 4. Applications

Consider the equation

$$z^{(r)} = f(z,t) + p(t), \quad r \geq 1, \tag{4.1}$$

$$f(z,t) = f(z,t+T), \quad f(0,t) = 0, \quad p(t) = p(t+T).$$

Assume that

$$f(z,t) = f(-z,-t), \quad p(t) = p(-t) \tag{4.2}$$

or

$$f(z,t) = f\left(-z, t + \frac{T}{2}\right), \quad p(t) = -p(t+T/2). \tag{4.3}$$

Under condition (4.2) or (4.3), equation (4.1) admits, respectively, a solution

$$z(t) = -z(-t) = z(t+T) \tag{4.4}$$

or

$$z(t) = -z\left(t + \frac{T}{2}\right). \tag{4.5}$$

These solutions have zero mean values and, therefore, are oscillatory.

According to Theorem 7, in the unperturbed system ($p(t) = 0$) under condition (2.34), $T$-periodic oscillatory solutions do not exist. The following theorem shows that for $p(t) \neq 0$, condition (2.34) guarantees the uniqueness of a solution (4.4) or (4.5).

**Theorem 8**. Under condition (2.34), systems (4.1), (4.2) and (4.1), (4.3) have no more than one solution of the kind (4.4) and (4.5), respectively.

**Proof**. Suppose that there are two such solutions, $z_1(t)$ and $z_{2_o}(t)$. The function $\delta(t) = z_2(t) - z_1(t)$ satisfies

$$\dot{\delta}(t) = C(t)\delta(t),  \qquad (4.6)$$

$$C(t) = \int_0^1 f_z(z_1(t) + s[z_2(t) - z_1(t)]), t)ds.$$

Obviously, the solution $\delta(t)$ is also of the kind (4.4) or (4.5), so, it is oscillatory. However, since $\|C(t)\| \leq N_L$, the existence of such solution in system (4.6) contradicts Theorem 7. □

Consider now the real linear Hamiltonian system with periodic coefficients,

$$\dot{x} = JH(t)x, \quad x \in \mathrm{R}^{2n}, \qquad (4.7)$$

$$J = \begin{bmatrix} 0 & -I \\ I & 0 \end{bmatrix}, \quad H(t) = H(t+T) = H^\tau(t) > 0.$$

Equation (4.7) is called strongly stable if all its solutions are bounded at infinity and retain this property under small perturbations of the Hamiltonian $H(t)$ which do not violate its symmetry [8].

In accordance with the Krein's theorem [1], for strong stability of equation (4.7) in the Hilbert space, it is sufficient that

$$\int_0^T \|JH(t)\| dt < \pi. \qquad (4.8)$$

Herewith, equation (4.7) belongs to the central stability region (all the multipliers of the first and second kind lie, respectively, on the upper and lower semicircles of the unit circle).

The following theorem generalizes this result to the space $\mathrm{R}^{2n}$ with any norm.

**Theorem 9**. System (4.7), (4.8) is strongly stable.

**Proof**. As is known [8], equation (4.7) belongs to the central stability region if $\lambda_1 > 1$, where $\lambda_1$ is the smallest positive eigenvalue of the boundary value problem

$$\dot{x} = \lambda JH(t)x, \quad x(T) = -x(0). \qquad (4.9)$$

At $\lambda = \lambda_1$, equation (4.9) has a multiplier $\rho = -1$ and, therefore, a $2T$-periodic solution $x_*(t) = -x_*(t+T)$. By (4.8) and $H(t) = H(t+T)$, for $\lambda \in (0,1]$,

$$\lambda \int_0^{2T} \|JH(t)\| dt < 2\pi.$$

As follows from Theorem 3, under this inequality, equation (4.9) has no a $2T$-periodic oscillatory solution. Consequently, $\lambda_1 > 1$, i.e., equation (4.7) is strongly stable. □

Consider the vector second order equation

$$\ddot{x} + P(t)x = 0, \quad x \in \mathbb{R}^n, \qquad P(t) > 0 \tag{4.10}$$

**Theorem 10**. Equation (4.10) is strongly stable, if for some $c > 0$,

$$\int_0^T \max[c, c^{-1}\|P(t)\|] dt < 2\pi. \tag{4.11}$$

**Proof.** As is known [8], a change of the variables reduces (4.10) to (4.7) with

$$H(t) = \begin{bmatrix} I & 0 \\ 0 & P(t) \end{bmatrix}.$$

Therefore, analogously to Theorem 9, the stability of equation (4.10) is guaranteed by the inequality $\lambda_1 > 1$. The last, according to Theorem 6, is provided by condition (2.26), which in the considered case, $r = 2$, takes the form (4.11). □

Compare (4.11) with the known stability condition for system (4.10) (Krein,,[1]) :

$$TR < 4, \tag{4.12}$$

$$R = \left\| \int_0^T P^+(t) dt \right\|, \quad P^+(t) = [|p_{ik}(t)|]_1^n.$$

Put $P(t, \varepsilon) = P_0 + \varepsilon P_1(t)$ where $\varepsilon > 0$ is a parameter. For $\varepsilon = 0$, the Euclidean norm, $\|P(t,0)\| = \|P_0\| = \lambda_0$ where $\lambda_0$ is the largest eigenvalue of the matrix $P_0$. Here the best bound (4.11) (reached for $c = \sqrt{\lambda_0}$) becomes.

$$T < \frac{\pi}{\sqrt{\lambda_0}} \tag{4.13}$$

This bound is sharp, because for $T = \pi/\sqrt{\lambda_0}$, equation (4.10) has a solution $x(t) = -x(t + T/2)$ and, therefore, is not strongly stable.

In turn, condition (4.12) for $\varepsilon = 0$ implies

$$T < \frac{2}{\sqrt{\lambda_+}} \leq \frac{2}{\sqrt{\lambda_0}} \qquad (4.14)$$

where $\lambda_+$ is the largest eigenvalue of the matrix $P^+(0)$.

As follows from (4.13) and (4.14), at least for small $\varepsilon$, the obtained stability condition (4.11) is less conservative than (4.12).

Consider the nonlinear Hamiltonian system

$$\dot{x} = JH_x(x,t), \quad x \in \mathrm{R}^{2n}, \qquad (4.15)$$

$$H_x(x,t) = H_x(x,t+T), \quad H_{xx}(x,t) > 0.$$

A solution $x(t) = x(t+T)$ of (4.15) is called stable to the first approximation, if the corresponding variational equation,

$$\dot{u} = JA(t)u, \quad A(t) = H_{xx}(x(t),t), \qquad (4.16)$$

is strongly stable.

Put

$$L_H(t) = \sup_{x \in \mathrm{R}^{2n}} \|JH_{xx}(x,t)\|, \quad N_H(T) = \int_0^T L_H(t)dt. \qquad (4.17)$$

**Theorem 11**. If

$$N_{HH}(T) < \pi, \qquad (4.18)$$

then the solution $x_1(t)$ is stable to the first approximation.

**Proof**. Obviously, $\|A(t)\| \leq L_H(t)$, therefore, under condition (4.18), the strong stability of equation (4.16) follows from Theorem 9. □

Note that if $H_x(0,t) = 0$, equation (4.15) has a solution $x(t) = 0$. So, here condition (4.18) guarantees stability of this equilibrium position.

Since the considered matrices $H(t)$ and $H_{xx}(x,t)$ are symmetric and $J^T J = I$, the matrices $JH(t)$ and $JH_{xx}(x,t)$ are normal. So, according to Lemma 1, the corresponding Euclidean norms have smallest values. Thus, in specific Hamiltonian systems, the best estimates are obtained when using the Euclidean norm.